\documentclass[reqno,12pt]{amsart}
\usepackage{amsmath, amssymb, amsthm,amsfonts} 
\usepackage[english]{babel}
\usepackage{bbm}
\usepackage{graphicx}
\usepackage{url}
\usepackage{epstopdf}\usepackage[a4paper,bindingoffset=0.5cm,left=2cm,right=2cm,top=2.5cm,bottom=2cm,footskip=.8cm]{geometry}

\usepackage{tikz}
\usetikzlibrary{arrows, automata,positioning,calc,shapes,decorations.pathreplacing,decorations.markings,shapes.misc,petri,topaths}
\usepackage{tkz-berge}
  \usepackage{pgfplots}
  \pgfplotsset{compat=newest}
  \usetikzlibrary{plotmarks}
  \usepackage{grffile}
\newlength\figureheight
  \newlength\figurewidth
\pgfplotsset{%
    tick label style={font=\scriptsize},
    label style={font=\footnotesize},
    legend style={font=\footnotesize},
         every axis plot/.append style={very thick}
}

\usepackage{rotating}
\usepackage{amsbsy,enumerate}
\usepackage{graphicx}
\usepackage{ccaption}
\usepackage{comment}
\usepackage{mathrsfs}

\newcommand{\s}{^\star}

\newcommand{\bs}{\boldsymbol}

\newcommand{\vb}{\vspace{3.2mm}}

\newcommand{\hr}{\hspace{5mm}}

\DeclareMathOperator*{\argmax}{arg\,max}

\allowdisplaybreaks

\newtheorem{corollary}{Corollary}

\newtheorem{theorem}{Theorem}
\newtheorem{remark}{Remark}

\newtheorem{algorithm}{Algorithm}

\makeatletter
\renewcommand{\fnum@figure}[1]{\textbf{\figurename~\thefigure}. }
\renewcommand{\fnum@table}[1]{\textbf{\tablename~\thetable}. }
\makeatother

\begin{document}

\title[A ruin model with resampling]{A ruin model with a resampled environment}

\author{C. Constantinescu, G. Delsing, M. Mandjes, L. Rojas Nandayapa}

\begin{abstract}
This paper considers a Cram\'er-Lundberg risk setting, where the components of the underlying model change over time. These components could be thought of as the claim arrival rate, the claim-size distribution, and the premium rate, but we allow the more general setting of the cumulative claim process being modelled as a spectrally positive L\'evy process. We provide an intuitively appealing mechanism to create such parameter uncertainty: at Poisson epochs we {\it resample} the model components from a finite number of $d$ settings. It results in a setup that is particularly suited to describe situations in which the risk reserve dynamics are affected by external processes (such as the state of the economy, political developments, weather or climate conditions, and policy regulations). 

\noindent We extend the classical Cram\'er-Lundberg approximation (asymptotically characterizing the all-time ruin probability in a light-tailed setting) to this more general setup. In addition, for the situation that the driving L\'evy processes are sums of Brownian motions and compound Poisson processes, we find an explicit uniform bound  on the ruin probability, which can be viewed as an extension of Lundberg's inequality; importantly, here it is not required that the L\'evy processes be spectrally one-sided. In passing we propose an importance-sampling algorithm facilitating efficient estimation, and prove it has bounded relative error. In a series of numerical experiments we assess the accuracy of the asymptotics and bounds, and illustrate that neglecting the resampling can lead to substantial underestimation of the risk.

\vb

\noindent
{\sc Keywords.} L\'evy risk processes $\circ$ parameter uncertainty $\circ$ ruin probabilities $\circ$ Cram\'er-Lundberg asymptotics $\circ$ Lundberg's inequality

\vb

\noindent
{\sc Affiliations.} 
{\it Corina Constantinescu} and {\it Leonardo Rojas Nandayapa} are with the Institute for Financial and Actuarial Mathematics, Department of Mathematical Sciences, University of Liverpool, L69 3 BX Liverpool, United Kingdom. 

Email: \url{{C.Constantinescu|leorojas}@liverpool.ac.uk}.

\noindent
{\it Guusje Delsing} is with Korteweg-de Vries Institute for Mathematics, University of Amsterdam, Science Park 904, 1098 XH Amsterdam, the Netherlands, and with Rabobank, Croeselaan 18, 3521 CB Utrecht, the Netherlands. Email: \url{g.a.delsing@uva.nl}.

\noindent
{\it Michel Mandjes} is with Korteweg-de Vries Institute for Mathematics, University of Amsterdam, Science Park 904, 1098 XH Amsterdam, the Netherlands. He is also with E{\sc urandom}, Eindhoven University of Technology, Eindhoven, the Netherlands, and Amsterdam Business School, Faculty of Economics and Business, University of Amsterdam, Amsterdam, the Netherlands. His research is partly funded by NWO Gravitation project N{\sc etworks}, grant number 024.002.003. Email: \url{m.r.h.mandjes@uva.nl}.

\end{abstract}

\maketitle

\newpage

\section{Introduction}
Risk theory focuses on analyzing models that describe an insurer's vulnerability to ruin.
Starting from the seminal works by Cram\'er \cite{C1} and Lundberg \cite{L1,L2} a substantial research effort has been spent on determining the ruin probability in a broad range of risk models. In the basic model, independent and identically distributed claims are assumed to arrive according to a Poisson process, whereas premiums arrive at a constant rate. The ruin probability is the probability that the capital surplus drops below $0$.

After the above mentioned pioneering papers, various extensions and generalizations have been considered to make the model more realistic. In this respect, multiple directions can be distinguished. Without pursuing to provide a complete overview, we include a brief account of a few important branches. In the first place, the classical model has been extended to include time-dependent ruin, i.e.\ ruin before a specified point in time; see e.g.\ \cite[Ch.\ V]{AA}. Secondly, the assumption of the cumulative claim process being of compound Poisson type has been generalized to that of compound Poisson perturbed by diffusion \cite{G2, G1}, and later to that of (spectrally one-sided) L\'evy input;  see e.g.\ \cite[Ch.\ X and XI]{AA} and \cite{KYP}. Thirdly, returns on investment have been included, and also level-dependent risk models have been considered; see e.g.\ \cite[Ch.\ VIII]{AA} and \cite{AC}. A major other branch in the literature focuses on computing or approximating ruin probabilities for specific claim-size distributions; see for instance \cite{CSZ} for the case of Gamma claims and \cite{RAM} for the case of heavy-tailed claims. Finally, we mention the direction of research in which the effect of specific dependence structures is assessed; see e.g.\ \cite{CKM} and, for an overview, \cite[Ch.\ XIII]{AA}. We also refer to  \cite{MIK,KGS, TEU} for further background on risk theory in general.

More often than not, in the models that have been considered the corresponding model primitives (in terms of parameters and distributions) are fixed. For instance in the classical Cram\'er-Lundberg model a specific claim arrival rate, premium rate, and claim-size distribution are held constant, in the sense that they cannot change over time. In reality, however, such a setup is typically not valid: as a consequence of various `external circumstances' the model primitives may fluctuate. In this context one could think of exogenous factors affecting the claim arrival process, such as the state of the economy, the political situation, weather conditions, and policy regulations.  Neglecting the parameter uncertainty (by using the conventional Cram\'er-Lundberg model with time-averaged parameters) could evidently lead to a substantial underestimation of the risk. 

An intuitively appealing mechanism to introduce parameter uncertainty is to periodically resample them. A very basic example of such a model would be an adaptation of the classical Cram\'er-Lundberg framework, in which (say every day, week or month) the arrival rate is resampled from a given distribution. Evidently, in principle also the other  model primitives (i.e., premium rate and claim-size distribution) can be periodically resampled. In  \cite{HLM}, for a different class of models, a similar mechanism to introduce parameter fluctuations has been proposed. 

In this paper we consider the setup in which the claim arrival process is a spectrally one-sided L\'evy process, thus covering the frequently used compound Poisson case. The special feature concerns the resampling mechanism described above: after exponentially distributed times, the Laplace exponent of this driving L\'evy process is resampled from a set of $d\in{\mathbb N}$ possible settings. There is a connection between this model and the one in which the claim arrival process is a so-called {\it Markov additive process} (MAP) \cite[Ch. VII]{AA}. Importantly, due to our specific resampling mechanism that we impose in this paper, the results we obtain are relatively explicit (compared to their counterparts under a MAP claim arrival process). Throughout this paper we assume the claim-size distributions are light-tailed (in line with what is assumed in the classical Cram\'er-Lundberg framework). 

The main contributions of our paper are the following. (i) In the first place, for an initial capital reserve level $u$, we identify the exact asymptotics of the ruin probability in the regime that $u$ grows large. This result can be seen as the counterpart of the Cram\'er-Lundberg asymptotics for our resampling model. (ii)~In the second place, restricting ourselves to the situation that the driving L\'evy processes are sums of Brownian motions and compound Poisson processes, we find an explicit upper bound  on the ruin probability that is uniform in $u\geqslant 0.$ This bound can be seen as an extension of Lundberg's inequality. In this context it is important to note that it is not required that the L\'evy processes be spectrally one-sided. (iii) In passing we propose an importance-sampling algorithm that facilitates the efficient estimation of small ruin probabilities. We prove that this procedure has bounded relative error, which effectively means that the number of runs needed to obtain an estimate of a given precision is hardly affected by the value of $u$. (iv) We conclude this paper by a series of numerical experiments, in which we systematically assess the accuracy of the asymptotics and bounds. An important observation is that neglecting the resampling (by using the Cram\'er-Lundberg model with time-averaged parameters) typically leads to a  significant underestimation of the risk.

This paper is organized as follows. Section \ref{PM} provides a formal model description and some preliminaries. Then in Section \ref{ANA} the exact asymptotics are established. Section \ref{UNI} presents the counterpart of Lundberg's inequality, together with the importance-sampling algorithm. Numerical examples are provided in Section \ref{EN}; this section also provides explicit expressions for the asymptotics and bounds in case the number of environmental states $d$ equals $2$. 

\section{Model and preliminaries}\label{PM}
In this section we introduce our resampling model, and provide preliminaries. In our model, the risk process is expressed in terms of a spectrally-positive L\'evy process, whose characteristics are resampled at Poisson epochs. 
\subsection{Model}
We start by constructing the net cumulative claim process $X(\cdot)$. To this end, 
we first introduce spectrally-positive scalar-valued L\'evy processes $X_i(\cdot)$ for $i=1,\ldots,d$, where we assume that $X_i(0)=0$. These processes are characterized by their respective Laplace exponents $\varphi_1(\cdot)$ up to $\varphi_d(\cdot)$, meaning that, for $\alpha\geqslant 0$,
\[\log {\mathbb E}\exp(-\alpha X_i(1)) = \varphi_i(\alpha);\]
see e.g.\ \cite{KYP}. 

In a standard ruin-theoretic setting the processes $X_i(\cdot)$ would correspond to compound Poisson processes (representing the  cumulative claim process) from which a deterministic drift is subtracted (the incoming premiums). Observe however that the framework we consider is significantly richer: the processes $X_i(\cdot)$ could contain a Brownian component, and also increasing `small-jumps processes' (such as the Gamma process or the Inverse Gaussian process) can be included \cite[Sections 1.2.4--1.2.5]{KYP}.

We now construct our resampling model. 
Let $T_n$ be the jump epochs of a Poisson process with rate $q>0$; we set $T_0:=0.$ 
At these epochs with probability $p_i\in[0,1]$ the $i$-th of the above-mentioned $d$ L\'evy processes is picked, with the $p_i$ summing to 1. Let $J_n\in\{1,\ldots,d\}$ be the index of the L\'evy process that was picked between $T_n$ and $T_{n+1}$, and set $J(t) = J_n$ when $t\in [T_n,T_{n+1})$.
Then we recursively define the cumulative claim process by, for $t\in [T_n,T_{n+1})$,
\[X(t) := X(T_n) + \big(X_{J_n}(t) - X_{J_{n}}(T_n)\big).\]
In a ruin context, we let $u-X(t)$ represent the capital surplus at time $t$, given the initial reserve was $u>0$. This means that the all-time ruin probability can be expressed as the probability that $X(t)\geqslant u$ for some $t\geqslant 0.$ This is the probability that we will study in this paper.

Define the all-time maxima
\[\bar X:=\sup_{t\geqslant 0} X(t),\:\:\: \bar X_{\rm d} := \sup_{n\in {\mathbb N}_0} X(T_n);\]
in other words, $X_{\rm d}$ is the all-time maximum, but restricted to jump epochs of the background process.
We work in the sequel with
\begin{align*}\pi(u):= &\:{\mathbb P}(\bar X\geqslant u)={\mathbb P}(\exists t\geqslant 0: X(t)\geqslant u),\:\:\:\\
\pi_{\rm d}(u):=&\: {\mathbb P}(\bar X_{\rm d}\geqslant u)=
{\mathbb P}(\exists n\in{\mathbb N}_0: X(T_n)\geqslant u).\end{align*}
It is clear that $\bar X\geqslant \bar X_{\rm d}$, so that $\pi(u)\geqslant \pi_{\rm d}(u).$ 

Throughout this paper we assume a negative drift, so that the events under consideration are increasingly rare as $u$ grows large. This negative drift assumption entails that we require
\begin{equation}\label{drift}\kappa:=\sum_{i=1}^d p_i\varphi'_i(0) >0.\end{equation}
In addition, in this work we assume that we are in the light-tailed setting, meaning that for all $i\in\{1,\ldots,d\}$ the Laplace exponent $\varphi_i(\alpha)$ is finite for $\alpha$ in an open neighborhood of the origin. In the $d=1$ case, this is in line with what was assumed to obtain the traditional Cram\'er-Lundberg asymptotics.

The claim arrival processes $X(\cdot)$  covers a resampled compound Poisson process as a special case.  
Then we can write the Laplace exponent of the $i$-th L\'evy process (i.e., $X_i(\cdot)$) as
\[\varphi_i(\alpha) = r_i\alpha -\lambda_i + \lambda_ib_i(\alpha),\]
where $r_i$ is the deterministic drift, $\lambda_i$ the claim arrival rate, and $b_i(\cdot)$ the Laplace transform of the claim sizes.

\subsection{Preliminaries}
In this paper the  focus lies in particular on the above probabilities' exact asymptotics (and related upper bounds) in the light-tailed domain. 
It is not hard to guess what the decay rate of the tail is. In the first place, one would expect that the logarithmic asymptotics of $\pi(u)$ and $\pi_{\rm d}(u)$ match (this we later prove). Secondly, observe that $(X(T_n))_{n\in {\mathbb N}}$ is a random walk; the increments $Y_n:=X(T_n)-X(T_{n-1})$ (for $n\in {\mathbb N}$) are independent and identically distributed (say, as a generic random variable $Y$). For this setting it is well-known \cite{KOR} that
\[\lim_{u\to\infty}\frac{1}{u}\log \pi_{\rm d}(u) = -\omega\s,\]
with $\omega\s$ the unique positive root of ${\mathbb E} \,{\rm e}^{\omega Y} = 1$. This means that $\omega\s$ solves
\begin{equation}
\label{omega}\sum_{i=1}^d p_i\int_0^\infty q{\rm e}^{-qt}{\rm e}^{\varphi_i(-\omega)\,t}{\rm d}t =
 \sum_{i=1}^d p_i\frac{q}{q-\varphi_i(-\omega)}=1
 \end{equation} 
 (where it is implicit that $\omega\s$ is such that $q>\varphi_i(-\omega\s)$ for all $i\in\{1,\ldots,d\}$). The existence of the root $\omega\s$ is assumed; it implies that there are $\alpha<0$ such that $\varphi_i(-\alpha)$  is finite (for all $i\in\{1,\ldots,d\}$), which means that we are in the regime that the upward jumps of $X_i(\cdot)$ are {\it light-tailed}; cf.\ e.g. \cite[Section 8.1]{DM}.
 
 Actually, the precise asymptotics of $\pi_{\rm d}(u)$ have been identified already. Recalling that $X(T_n)$ can be written as the sum of $n$ independent and identically distributed increments $Y_1$ up to $Y_n$, the exceedance probability $\pi_{\rm d}(u)$ can be interpreted as the probability that a random walk with negative drift (cf.\ condition (\ref{drift})) and light-tailed increments ever exceeds level $u$. For this setting in e.g.\ \cite{KOR} a positive constant $\gamma$ is found  that $\pi_{\rm d}(u)\,{\rm e}^{\omega\s u} \to \gamma$. As these exact asymptotics of $\pi(u)$ have not been identified so far, it is one of the main objectives of this paper to derive these; see Section \ref{ANA}. Another objective concerns a uniform upper bound on $\pi(u)$; see Section~\ref{UNI}.

\section{Asymptotics}\label{ANA}
In order to identify the  exact asymptotics of $\pi(u)$, we first verify that our model actually corresponds to the maximum value attained by a specifically chosen Markov additive process (in the sequel abbreviated to MAP); see e.g.\  \cite[Section 11.4]{DM}. To this end, recall that a MAP behaves as a L\'evy process $X_i(\cdot)$ whenever the background process $J(\cdot)$ (whose  transition rate matrix we denote by $Q=(q_{ij})_{i,j=1}^d$) is in state $i\in\{1,\ldots,d\}$. Let us construct the matrix $Q$, by considering the transition rates from state $i$. Observe that the process  $J(\cdot)$ stays for an exponential amount of time (with rate, say, $\bar q_i$) in $i$;  after this time, it jumps to  state $j\not=i$ with probability $p_j/(1-p_i)$. The parameter  $\bar q_i$ can be determined by computing the Laplace-Stieltjes transform of the time spent in state $i$, say $\tau_i$. We obtain
\[{\mathbb E}\,{\rm e}^{-\alpha\tau_i}=\sum_{k=1}^\infty p_i^{\,k-1}\left(1-p_i\right)\left(\frac{q}{q+\alpha}\right)^k
=\frac{(1-p_i)\,q}{\alpha+(1-p_i)\,q},\]
from which we conclude that $\tau_i$ is exponentially distributed with parameter $\bar q_i=(1-p_i)\,q$. We thus observe that 
$q_{ij} = \bar q_j \cdot p_j/(1-p_i) = q\,p_j$  for $i\not=j$, whereas $q_{ii}= -\bar q_i$. We thus arrive at  
\begin{equation}\label{Q}Q = q \,{\bs e}{\bs p}^\top -qI_d,\end{equation} with ${\bs e}$ an all-ones vector and $I_d$ the $d$-dimensional identity matrix. The conclusion is that our process $X(\cdot)$ corresponds to a MAP with the transition rate matrix $Q$  given by (\ref{Q}), and state-dependent Laplace exponents $\varphi_i(\cdot)$ for $i\in\{1,\ldots,d\}.$ In the sequel, we use that $X(\cdot)$ has a MAP-representation; in particular we make use of the fact that for this setting the Laplace transform of $\bar X$ is known.

\subsection{Transform of $\bar X$} In this subsection we provide the Laplace transform of $\bar X$, and show how this can be simplified, owing to the special structure of the matrix $Q$.

Define $\Phi(\alpha):={\rm diag}\{\varphi_1(\alpha),\dots,\varphi_d(\alpha)\}$ and  $\bar\Phi(\alpha) =-qI_d + \Phi(\alpha)$. We can now introduce 
\[M(\alpha) = Q + \Phi(\alpha) = q \,{\bs e}{\bs p}^\top-qI_d + \Phi(\alpha) = q \,{\bs e}{\bs p}^\top+ \bar\Phi(\alpha),\]
which can be seen as a `matrix-valued Laplace exponent' in the sense that, for any $i,j\in\{1,\ldots,d\}$,
\[{\mathbb E}\big({\rm e}^{-\alpha X(t)}1_{\{J(t)=j\}}\,|\,J(0)=i\big) = ({\rm e}^{M(\alpha) t})_{i,j}.\]
By a Perron-Frobenius based argumentation, one can show that the eigenvalue of $M(\alpha)$ with largest real part, which we denote by $\mu(M(\alpha))$, is actually real (where we note that in our specific setting we  argue below that {\it all} $d$ eigenvalues are real). It thus follows that
\[\lim_{t\to\infty}\frac{1}{t}\log {\mathbb E}{\rm e}^{\alpha X(t)} = \mu(M(\alpha)).\]
Due to the fact that this concerns a limiting logarithmic moment generating function, we thus conclude that $\mu(M(\alpha))$ is a convex function of $\alpha$.

The Laplace transform of $\bar X$ in $\alpha$ can also be expressed in terms of this matrix $M(\alpha).$ More specifically, 
as can be found in e.g.\ \cite[Eqn.\ (11.1)]{DM}, there is the following `matrix counterpart' of the celebrated Pollaczek-Khinchine formula:
\[{\mathbb E}\,{\rm e}^{-\alpha \bar X} = \alpha\,{\bs\ell}^\top \,(M(\alpha))^{-1}{\bs e},\]
for a  vector ${\bs\ell}$ determined in e.g.\ \cite{DIKM,DiM}. 

\begin{remark}\label{R1}{\em In \cite{DIKM} a compact representation for the vector ${\bs\ell}$ is given. We provide a brief account of this representation here. First, split the $d$ background states as follows. Let for states $i\in \{d-d^-+1,\ldots,d\}$ the L\'evy process $X_i(\cdot)$  correspond to a decreasing subordinator; obviously, in these states the process $X(\cdot)$ cannot cross the level $u$ (if there are no states corresponding to decreasing subordinators, we put $d^-:=0$). In the other states, corresponding to $i\in\{1,\ldots,d-d^-\}$, the level $u$ {\it can} be crossed. 

Now consider
\[\eta(v):=\inf\{t\geqslant 0: -X(t) \geqslant v\},\]
as a process in $v\geqslant 0.$
As argued in e.g.\ \cite{DIKM}, $(\eta(v), J(\eta(v))_{v\geqslant 0}$ is a MAP, with $J(\eta(v))$ attaining values in $\{1,\ldots,d-d^-\}$. Let $\bar{\bs\pi}$ be the $(d-d^-)$-dimensional invariant probability measure pertaining to the Markov process 
$(J(\eta(v))_{v\geqslant 0}$.
Then the vector ${\bs\ell}$ is such that ${\bs\ell}^\top= \kappa\,(\bar{\bs\pi}^\top,{\bs 0}^\top)$, with the scalar $\kappa>0$ as in (\ref{drift}).\hfill$\Box$
}\end{remark}

Interestingly, due to the fact that $M(\alpha)$ is the sum of a diagonal matrix and a rank-one matrix, its eigenvalues can be somehow characterized, applying the following nice (and well-known) idea. To this end, we can write 
\[{\rm det}(M(\alpha)-\theta I_d) = {\rm det}(\bar \Phi(\alpha)-\theta I_d){\rm det}(I_d+(\bar \Phi(\alpha)-\theta I_d)^{-1} q \,{\bs e}{\bs p}^\top).\]
For $A$ of dimension $m\times n$, and $B$ of dimension $n\times m$, we have ${\rm det}(I_m-AB)={\rm det}(I_n-BA)$. We thus conclude that
\begin{align*}
{\rm det}(M(\alpha)-\theta I_d) =&\: {\rm det}(\bar \Phi(\alpha)-\theta I_d)\,{\rm det}(I_d+\,{\bs p}^\top(\bar \Phi(\alpha)-\theta I_d)^{-1} q \,{\bs e})\\
=&\: {\rm det}(\bar \Phi(\alpha)-\theta I_d)\left(1-p_i\sum_{i=1}^d \frac{q}{q-\varphi_i(\alpha)+\theta}\right).
\end{align*}
We find that the eigenvalues  $\theta_1(\alpha)$ up to $\theta_d(\alpha)$, for a fixed $\alpha$, are the solutions to
\[\frac{1}{q} = \sum_{i=1}^d p_i \frac{1}{q-\varphi_i(\alpha)+\theta}=:\Psi_\alpha(\theta).\]
With $\Theta(\alpha):={\rm diag}\{\theta_1(\alpha),\ldots,\theta_d(\alpha)\}$, using standard machinery from linear algebra we get for a matrix $S(\alpha)$ that
\[{\mathbb E}\,{\rm e}^{-\alpha \bar X} = \alpha{\bs\ell}^\top S(\alpha) (\Theta(\alpha))^{-1} S^{-1}(\alpha)\,{\bs e},\]
under the familiar regularity conditions regarding the multiplicities of the eigenvalues. In principle we have now a unique characterization of $\bar X$, and hence also, albeit in implicit terms, a way of computing $\pi(u).$ In general this requires numerical inversion, for which there are various algorithms available; see e.g.\ \cite{AW,dI}. In this section we have another objective: we use knowledge of the transform of $\bar X$ to identify the corresponding tail asymptotics. 

\begin{remark}\label{R2}
{\em A standard fact from linear algebra is that the columns of $S(\alpha)$ contain the right eigenvectors of $M(\alpha).$ If the eigenvalues $\theta_1(\alpha)$ up to $\theta_d(\alpha)$ have been found, these can be easily expressed in terms of these eigenvalues. Suppose $\theta$ is such an eigenvalue. Then the eigenvector ${\bs x}$ satisfies $M(\alpha){\bs x} = \theta {\bs x},$ or, equivalently, for $j\in\{1,\ldots,d\}$,
\[q\,{\bs p}^\top{\bs x} - qx_j +\varphi_j(\alpha)\,x_j =\theta x_j.\]
We conclude that
\[\frac{x_j}{x_i} =\frac{q-\varphi_i(\alpha)+\theta}{q-\varphi_j(\alpha)+\theta},\]
so that we can pick $x_i= (q-\varphi_i(\alpha)+\theta)^{-1}$.\hfill$\Box$
}
\end{remark}

\subsection{Tail asymptotics} The idea is to rely on the Heaviside recipe \cite[Recipe 8.1]{DM} to find the tail behavior. To this end we first have to identify the rightmost pole on the negative halfline. The poles are the values of $\alpha<0$ for which one of the $\theta_i(\alpha)$ equals 0. 

\begin{figure}
\begin{tikzpicture}
      \draw[->] (-3,0) -- (2.2,0) node[right] {$\theta$};
      \draw[->] (0,-2.3) -- (0,2.3) node[above] {$\Psi_\alpha(\theta)$};
      \draw[scale=0.5,domain=-5:-2.08,smooth,variable=\x,blue] plot ({\x},{1/3*1/(2+\x)
      +1/3*1/(1+\x)+1/3*1/(\x-1)});
      \draw[scale=0.5,domain=-1.92:-1.08,smooth,variable=\x,blue] plot ({\x},{1/3*1/(2+\x)
      +1/3*1/(1+\x)+1/3*1/(\x-1)});
      \draw[scale=0.5,domain=-0.92:0.92,smooth,variable=\x,blue] plot ({\x},{1/3*1/(2+\x)
      +1/3*1/(1+\x)+1/3*1/(\x-1)});
      \draw[scale=0.5,domain=1.08:3,smooth,variable=\x,blue] plot ({\x},{1/3*1/(2+\x)
      +1/3*1/(1+\x)+1/3*1/(\x-1)});
 \draw[scale=0.5,domain=-4:4,smooth,variable=\y,red]  plot ({-2},{\y});
  \draw[scale=0.5,domain=-4:4,smooth,variable=\y,red]  plot ({-1},{\y});
   \draw[scale=0.5,domain=-4:4,smooth,variable=\y,red]  plot ({1},{\y});
   \draw[dotted][scale=0.5,domain=-5:4,smooth,variable=\y,black]  plot ({\y},{1});
    \end{tikzpicture}
    \caption{\label{F1}For a given value of $\alpha$, a plot of $\Psi_\alpha(\theta)$ as a function of $\theta$. It illustrates the statement on the locations of the eigenvalues. The red vertical lines correspond to the poles $\bar\varphi_i(\alpha) $. The dotted horizontal line is at level $1/q.$}
    \end{figure}
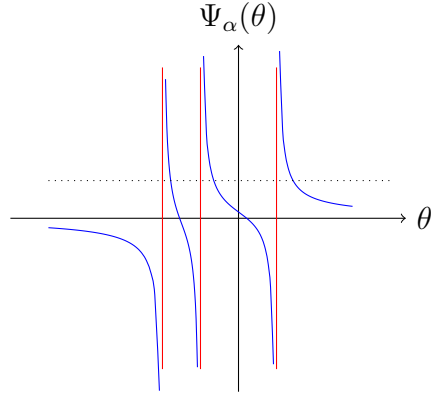

To study the behavior of the poles, first observe that $\Psi_0(0) = 1/q$, so for $\alpha=0$ all $d$ roots equal $0$. Now pick a negative value of $\alpha$, and let the bijection $b(i,\alpha)$ relabel the $\varphi_i(\alpha)$ such that, with 
\[\bar\varphi_i(\alpha):=\varphi_{b(i,\alpha)}(\alpha)-q,\]
the $\bar\varphi_i(\alpha)$ are increasing in $i$. Then, using the shape of $\Psi_\alpha(\theta)$, it is easily argued that one eigenvalue is larger than $\bar\varphi_d(\alpha)$, and that for $i=1,\ldots,d-1$ there is one of the eigenvalues $\theta_j(\alpha)$ in each of the intervals $(\bar\varphi_{i}(\alpha),\bar\varphi_{i+1}(\alpha))$, as illustrated in Fig.\ \ref{F1}; to this end, observe that
\[\lim_{\theta\uparrow \bar\varphi_i(\alpha) }\Psi_\alpha(\theta) = -\infty,\:\:\:
\lim_{\theta\downarrow \bar\varphi_i(\alpha) }\Psi_\alpha(\theta) = \infty,\]
and
\[\lim_{\theta\to-\infty}\Psi_\alpha(\theta) = \lim_{\theta\to\infty}\Psi_\alpha(\theta) = 0.\]
In this argumentation it is tacitly assumed that the $\bar\varphi_i(\alpha)$ are different, but the reasoning followed extends in an obvious way to the case that some are equal.

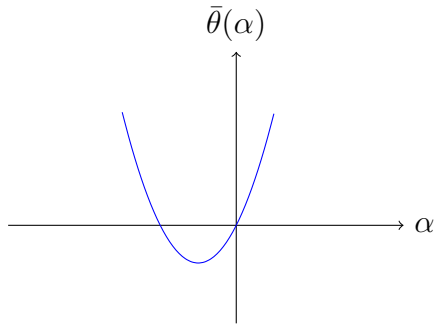
\begin{figure}
\begin{tikzpicture}
      \draw[->] (-3,0) -- (2.2,0) node[right] {$\alpha$};
      \draw[->] (0,-1.3) -- (0,2.3) node[above] {$\bar\theta(\alpha)$};
      \draw[scale=0.5,domain=-3:0.99,smooth,variable=\x,blue] plot ({\x},{\x*(2+\x)}); 
    \end{tikzpicture}
    \caption{\label{F2}The function $\bar\theta(\alpha)=\mu(M(\alpha))$.}
    \end{figure}

Denote $\bar\theta(\alpha):=\max_{i=1,\ldots,d}\theta_i(\alpha)$, which equals the $\mu(M(\alpha))$ we introduced earlier. As we observed above, $\bar\theta(0)=0.$
\begin{itemize}
\item[$\circ$] First consider the case of $\alpha=-\epsilon$ in the regime $\epsilon\downarrow 0.$ It is readily see that we are faced with the equation, putting $\theta=\delta\epsilon$ (so that $\delta=-\bar\theta'(0)$),
\[\frac{1}{q}= \sum_{i=1}^d p_i\frac{1}{q+\epsilon \varphi_i'(0)+\delta\epsilon+O(\epsilon^2)}=\frac{1}{q}\sum_{i=1}^d p_i\left(1-\frac{\varphi_i'(0)+\delta}{q}\epsilon\right)+O(\epsilon^2).\]
This leads to
\[\delta = -\sum_{i=1}^d p_i\varphi'(0),\]
which we know is negative due to the drift condition $(\ref{drift}).$
Conclude that $\bar\theta(\alpha)<0$ for small negative $\alpha$.
\item[$\circ$] We assume $\varphi_i(\alpha)\to\infty$ as $\alpha\to-\infty$ for at least one $i$ (to avoid trivial cases). As we know that there is one eigenvalue larger than $\bar\varphi_d(\alpha)$, we conclude that for $\alpha$ below some negative threshold, $\bar\theta(\alpha)>0$.
\item[$\circ$]
Recall that, from the interpretation of $\bar\theta(\alpha)$ as the limiting log moment generating function $\mu(M(\alpha))$, we know it is convex; see Figure \ref{F2}. 
 \end{itemize} As a consequence of these observations, we now conclude that the `rightmost pole on the negative halfline' is well defined, and characterized as 
\[\omega\s := -\sup\{\alpha<0: \bar\theta(\alpha)=0\},\]
which solves (\ref{omega}).

Now that we have identified the rightmost pole on the negative halfline, we are in a position to apply the Heaviside approach. To this end, we first note that
\[\zeta(\alpha):=\int_0^\infty e^{-\alpha u} \pi(u) {\rm d} u = \frac{1}{\alpha}\big(1-{\mathbb E}\,{\rm e}^{-\alpha\bar X}\big)
=\frac{1}{\alpha}-{\bs\ell}^\top S(\alpha) (\Theta(\alpha))^{-1} S^{-1}(\alpha)\,{\bs e}.\] We introduce \[i\s:=\argmax_{i=1,\ldots,d} \theta_i(-\omega\s),\]
which is the index of the eigenvalue that corresponds to the pole in $-\omega\s$.
Now define
\[A:=\lim_{\alpha\downarrow -\omega\s} (\alpha+\omega\s)\,\zeta(\alpha).\]
Then the Heaviside principle entails that
\[\pi(u){\rm e}^{\omega\s u}\to A = -{\bs\ell}^\top S(-\omega\s) \left(\lim_{\alpha\downarrow -\omega\s}(\alpha+\omega\s)(\Theta(\alpha))^{-1}\right) S^{-1}(-\omega\s)\,{\bs e}.\]
Denote by ${\bs u}$ the $i\s$-th column of $S(-\omega\s)$ and by ${\bs v}$ the $i\s$-th row of $S^{-1}(-\omega\s)$. Then
\begin{equation}
\label{AA}A = -\big({\bs\ell}^\top {\bs u}\big)\big( {\bs v}^\top{\bs e} \big)\frac{1}{\theta_{i\s}'(-\omega\s)}.\end{equation}
By Remark \ref{R2}, we have $u_i = (q-\varphi_i(-\omega\s))^{-1}$.

As pointed out in e.g.\ \cite[Section 3]{AW1},
the Heaviside recipe can be rigorously justified in some cases, but remains to be in others; we refer to e.g.\   Doetsch \cite[p.\ 254]{DOE}  and \cite[Sections 3 and 5]{ACW} for in-depth technical discussions. Importantly, for the case of the maximum of a spectrally-positive L\'evy process (without resampling, that is), it is argued in \cite[Section 8.1]{DM} that applying the Heaviside recipe to the generalized Pollaczek-Khinchine formula \cite{ZO} indeed yields the correct exact asymptotics; these asymptotics were derived in e.g.\ \cite{BD}, and can be seen as  the extension of the classical Cram\'er-Lundberg asymptotics to the case that the  L\'evy process is spectrally-positive. For our model we also assume that the use of the Heaviside recipe is justified. 

The above leads to the following generalization of the classical Cram\'er-Lundberg asymptotics.
\begin{theorem}\label{THCL}
As $u\to\infty$, $\pi(u) \,{\rm e}^{\omega\s u}\to A$, with $A$ given by $(\ref{AA})$.
\end{theorem}

\begin{corollary}
The probabilities $\pi(u)$ and $\pi_{\rm d}(u)$ are asymptotically proportional, in that their ratio tends to a positive constant as $u\to\infty.$
\end{corollary}

\begin{remark}\label{RA}{\em In the literature, results related to Thm.\  \ref{THCL} have appeared. We refer to e.g.\ \cite[Thm.\ 3.7]{AA} for a setting covering the L\'evy processes being compound Poisson processes.}
\end{remark}

\section{Uniform bound, change-of-measure, importance sampling} 
\label{UNI}
An intrinsic drawback of the asymptotics presented in the previous section is that they apply for large $u$ only; in addition, no explicit error bounds are provided. As a result, we do not know how accurate (for a given value of $u$) the approximation $\pi(u)\approx A\,{\rm e}^{-\omega\s u}$ is. This observation motivates the interest in searching for an upper bound that holds {\it uniformly in}~$u$. Here it is noted that, with the application in ruin theory, determining the initial capital level using an upper bound on $\pi(u)$ has the attractive feature that it leads to a `safe' policy. The main finding of this section is an upper bound on $\pi(u)$ which is proportional to ${\rm e}^{-\omega\s u}$, with the constant $\omega\s>0$ as defined before (i.e., as the solution of ${\mathbb E}\,{\rm e}^{\omega Y}=1$).
The bound can be seen as an extension of the classical Lundberg's inequality \cite[Thm. IV.5.2]{AA} to our model with resampling. 

Our proof, leading to the uniform upper bound in Thm.\ \ref{THUB}, is based on a change-of-measure argument. As a result, the reasoning also reveals in passing how importance sampling can be performed.  In Thm.\ \ref{BRE} we show that the importance procedure proposed is endowed with bounded relative error.

In this section, we consider the situation that the driving L\'evy processes are sums of Brownian motions and compound Poisson processes (i.e., do not include components with infinitely many `small jumps'); at the same time, we lift the assumption that the processes $X_i(\cdot)$ be spectrally positive. In practical terms, the fact that our L\'evy processes are not allowed to have a small jumps part is not a real restriction. As pointed out in \cite[Ch. X]{DM}, in simulation one could approximate the small jumps components by appropriately chosen Brownian motions, based on results in e.g.\ \cite{AR}.

Define $\tau(u):=\inf\{t\geqslant 0: X(t) \geqslant u\},$
such that $\pi(u)={\mathbb P}(\tau(u)< \infty).$
We proceed by analyzing this probability under a particular alternative measure ${\mathbb Q},$ defined as follows. The measure  ${\mathbb Q}$ is constructed such that (in self-evident notation), with $Y$ as defined before,
\[{\mathbb E}_{\mathbb Q} {\rm e}^{\omega Y}=\frac{{\mathbb E}\,{\rm e}^{(\omega+\omega\s) Y}}{{\mathbb E}\, {\rm e}^{\omega\s Y}}={\mathbb E} \, {\rm e}^{(\omega+\omega\s) Y},\]
where the second equality is by the definition of $\omega\s.$
Rewriting the right-hand side of the previous display as 
\[\sum_{i=1}^d p_i \frac{q}{q-\varphi_i(-\omega\s)}\left(\frac{q-\varphi_i(-\omega\s)}{q-\varphi_i(-\omega\s)
-\varphi_i(-\omega-\omega\s) +\varphi_i(-\omega\s)}\right),\]
and comparing with (\ref{omega}), 
we observe that we should choose, for $i\in\{1,\ldots,d\}$,
\begin{align*}
p_i ^{\mathbb Q}&:= p_i \frac{q}{q-\varphi_i(-\omega\s)},\:\:\:\:\:
q_i^{\mathbb Q}:= q-\varphi_i(-\omega\s),\:\:\:\:\:
\varphi_i^{\mathbb Q}(\cdot) := \varphi_i(\cdot\,-\omega\s) -\varphi_i(-\omega\s) ,
\end{align*}
with the superscript ${\mathbb Q}$ denoting that the parameters correspond to the new measure.

We now detail how the parameters of the Brownian motions and compound Poisson processes should be adapted under ${\mathbb Q}.$ We can write the Laplace exponent of the $i$-th L\'evy process (under the original measure) as
\[\varphi_i(\alpha) = r_i\alpha +\tfrac{1}{2}\sigma_i^2 \alpha^2 -\lambda_i + \lambda_ib_i(\alpha),\]
where $r_i$ is the deterministic drift, $\sigma_i^2$ the variance pertaining to the Brownian motion, $\lambda_i$ the claim arrival rate, and $b_i(\cdot)$ the Laplace transform of the claim sizes (where, as mentioned above, negative claims are allowed). 
We note that $\varphi_i^{\mathbb Q}(\cdot)$ is \cite[Section 10.2]{DM} the Laplace exponent of a L\'evy process (exponentially twisted with parameter $\omega\s$, that is); actually it is a sum of a Brownian motion and a compound Poisson process. It takes a minor computation to verify that $\varphi_i^{\mathbb Q}(\cdot) = \varphi_i(\cdot\,-\omega\s) -\varphi_i(-\omega\s)$ translates into (in self-eivident notation)
\[r_i^{\mathbb Q} = r_i-\omega\s\sigma_i^2,\:\:\:\:\lambda_i^{\mathbb Q}:=\lambda_ib_i(-\omega\s),\:\:\:\: b_i^{\mathbb Q}(\cdot):=\frac{b_i(\cdot\,-\omega\s)}{b_i(-\omega\s)},\]
where $\sigma_i^2$ remains unchanged.

Observe that (i)~due to the definition of $\omega\s$ the $p_i$\,s sum to $1$, (ii) $q_i^{\mathbb Q}>0$ for all $i\in\{1,\ldots,d\}$ (recall that $\omega\s$ is such that $q>\varphi_i(-\omega\s)$). Note that under ${\mathbb Q}$ the times spent in the states $1$ up to $d$ are still exponential, but now with a state-dependent parameter (whereas this parameter was state-independent under ${\mathbb P}$); informally, the measure ${\mathbb Q}$ increases the preference for states under which ruin is relatively likely.

Due to the convexity of moment generating functions,
\[{\mathbb E}_{\mathbb Q} Y = \left.\frac{\rm d}{{\rm d}\omega} {\mathbb E}_{\mathbb Q} {\rm e}^{\omega Y}\right|_{\omega=0}=\left.\frac{\rm d}{{\rm d}\omega} \frac{{\mathbb E}\, {\rm e}^{(\omega+\omega\s) Y}}{{\mathbb E}\, {\rm e}^{\omega\s Y}}\right|_{\omega=0}=\left.\frac{\rm d}{{\rm d}\omega} {{\mathbb E} \, {\rm e}^{(\omega+\omega\s) Y}}\right|_{\omega=0}>0,\]
which implies that ${\mathbb Q}(\tau(u)<\infty)=1.$ In other words, we have constructed a new measure under which the event under consideration happens almost surely. We thus find the identity \cite[Section XIII.3]{ASM}
\[\pi(u) ={\mathbb P}(\tau(u)<\infty) = {\mathbb E}_{\mathbb Q}L,\]
where $L$ is the likelihood ratio (under ${\mathbb P}$, relative to ${\mathbb Q}$, that is) corresponding to the trajectory of the stochastic process $X(\cdot)$ until $u$ has been reached (i.e., time $\tau(u)$). One could write
\[L=\frac{{\rm d}{\mathbb P}}{{\rm d}{\mathbb Q}}\big((X(t))_{t\in[0,\tau(u)]}\big).\]

The next observation is that $u$ is (first) reached either (i)~due to Brownian motion attaining the value $u$ in between two consecutive claim arrivals, or (ii)~due to a claim arrival. Supposing that  at some point in time the background state is $i$, the time till either a change of the background state or a claim arrival is exponential with parameter $f^{\mathbb Q}_i:= \lambda^{\mathbb Q}_i+q_i^{\mathbb Q}$ (which used to be $f_i:= \lambda_i+q$ under the original measure). The increment of the process $(X_t)_{t\geqslant 0}$ in this interval can be written as the sum of three independent terms:
\begin{itemize}
\item[$\circ$]
In the first place there is the maximum attained by the Brownian motion in the interval. This is a positive term, that is exponentially distributed (under ${\mathbb Q}$, that is)
with parameter
\[\alpha_{i,+}^{\mathbb Q}:= \frac{\sqrt{(r_i^{\mathbb Q})^2+2f_i^{\mathbb Q}\sigma_i^2}+r_i^{\mathbb Q}}{\sigma_i^2}.\]
\item[$\circ$] In the second place there is the (negative) distance between this maximum and the value at the end of the interval, just prior to the claim arrival. This is a negative term, of which the absolute value is exponentially distributed (under ${\mathbb Q}$) with parameter
\[\alpha_{i,-}^{\mathbb Q}:= \frac{\sqrt{(r_i^{\mathbb Q})^2+2f_i^{\mathbb Q}\sigma_i^2}-r_i^{\mathbb Q}}{\sigma_i^2}.\]
\item[$\circ$] In the third place there is the claim size, which is sampled from a distribution with Laplace transform $b_i^{\mathbb Q}(\cdot).$
\end{itemize}
The justification of the above decomposition (and, in particular, the independence between the first two terms) lies in Wiener-Hopf arguments; see e.g.\ \cite[Ch. VI]{KYP}. More specifically, we have the following expression for the Laplace transform of the value of the Brownian component of $X_i(\cdot)$ after an exponentially distributed interval with mean $f_i^{-1}$ under the original measure:
\begin{align*}\int_0^\infty f_i \,{\rm e}^{-f_it}\,{\mathbb E}\,{\rm e}^{-\alpha X_i(t)}{\rm d}t&=\frac{f_i}{f_i-\varphi_i(\alpha)}
=\frac{f_i}{f_i-r_i\alpha -\tfrac{1}{2}\sigma_i^2\alpha^2}\\
&=\frac{2f_i}{\sigma_i^2}\cdot\frac{1}{\alpha_{i,+}+\alpha}\cdot\frac{1}{\alpha_{i,-}-\alpha},
\end{align*}
with
\[\alpha_{i,+}:=\frac{\sqrt{r_i^2+2f_i\sigma_i^2}+r_i}{\sigma_i^2},\:\:\: \alpha_{i,-}:=\frac{\sqrt{r_i^2+2f_i\sigma_i^2}-r_i}{\sigma_i^2},\]
which are both positive numbers;  an analogous reasoning applies under ${\mathbb Q}$.

We now present (in self-evident notation) a pseudocode for the importance sampling procedure that we propose. We let $B_i$ represent i.i.d.\ samples from a distribution with Laplace transform $b^{\mathbb Q}_i(\cdot)$, $A^+_i$ are i.i.d.\ samples from an exponential distribution with mean $1/\alpha^{\mathbb Q}_{i,+}$, and $A^-_i$ are i.i.d.\ samples from an exponential distribution with mean $1/\alpha^{\mathbb Q}_{i,-}$.
As mentioned, termination (i.e., reaching `{\tt \small EXIT}') of the algorithm is guaranteed by ${\mathbb E}_{\mathbb Q}Y>0.$
Each time `{\tt\small Random}' appears in the algorithm, a new (i.e., independent of all previous ones) uniform random number is generated (on the interval $[0,1]$). 

\begin{algorithm} \label{ALG1} {\em \tt \small $X:=0$;\:\:\:$L:=1$; 

\noindent REPEAT

\hr $I \sim {\bs p}^{\mathbb Q}$;\:\:\:$L:=L * p_I/p_I^{\mathbb Q}$;

\hr WHILE Random $<\, \lambda^{\mathbb Q}_I/f^{\mathbb Q}(I)$ THEN

\hr\hr $L := L * (f^{\mathbb Q}_I/f_I) * (\lambda_I/\lambda_I^{\mathbb Q})$;

\hr\hr ${\mathbb A}^+\sim A^+_I$;\:\:\: $X := X+ {\mathbb A}^+$; \:\:\:$L:= L* (\alpha_{I,+}/\alpha_{I,+}^{\mathbb Q}) * \exp(-(\alpha_{I,+}-\alpha_{I,+}^{\mathbb Q}){\mathbb A}^+)$;

\hr\hr IF $X>u$ THEN RETURN $L$;\:\:\:EXIT;

\hr\hr ${\mathbb A}^-\sim A^-_I$;\:\:\: $X := X- {\mathbb A}^-$; \:\:\:$L:= L* (\alpha_{I,-}/\alpha_{I,-}^{\mathbb Q}) * \exp(-(\alpha_{I,-}-\alpha_{I,-}^{\mathbb Q}){\mathbb A}^-)$;

\hr\hr ${\mathbb B}\sim B_I$;\:\:\: $X:=X +{\mathbb B}$; \:\:\:$L:= L * \exp(-\omega\s {\mathbb B}) * b_I(-\omega\s)$;

\hr\hr IF $X>u$ THEN RETURN $L$;\:\:\:EXIT;

\hr END (of `WHILE');

\hr $L := L * (f^{\mathbb Q}_I/f_I) * (q_I/q_I^{\mathbb Q})$;

\hr ${\mathbb A}^+\sim A^+_I$;\:\:\: $X := X+ {\mathbb A}^+$; \:\:\:$L:= L* (\alpha_{I,+}/\alpha_{I,+}^{\mathbb Q}) * \exp(-(\alpha_{I,+}-\alpha_{I,+}^{\mathbb Q}){\mathbb A}^+)$;

\hr IF $X>u$ THEN RETURN $L$;\:\:\:EXIT;

\hr ${\mathbb A}^-\sim A^-_I$;\:\:\: $X := X- {\mathbb A}^-$; \:\:\:$L:= L* (\alpha_{I,-}/\alpha_{I,-}^{\mathbb Q}) * \exp(-(\alpha_{I,-}-\alpha_{I,-}^{\mathbb Q}){\mathbb A}^-)$;

UNTIL FALSE. \hfill$\Box$

}\end{algorithm}

Let us now evaluate $L$, as resulting from Algorithm \ref{ALG1}, in greater detail; we do these computations to derive a uniform upper bound on $\pi(u)$. 
Define the variable $N$ as the number of times the background state is resampled in the simulation until level $u$ is reached; equivalently, $u$ is reached in $(T_{N-1},T_N].$ Now consider the contribution $L_n$ to the likelihood ratio $L$ due to the random objects sampled in the interval $(T_{n-1},T_n]$, for $n\in\{1,\ldots,N\}$; as a consequence, $L=L_1\cdots L_N.$ 

We state by considering $n\in\{1,\ldots,N-1\}$. Let there have been $K_n$ claim arrivals in that interval; let $I_n$ be the background state in this interval. Let $(A_{i,j}^\pm)_{j\geqslant 1}$ be i.i.d.\ copies of $A_i^\pm$, and $(B_{i,j})_{j\geqslant 1}$ i.i.d.\ copies of $B_i$. Then, in self-evident notation,
\begin{align*} L_n&= \frac{p_{I_n}}{p_{I_n}^{\mathbb Q}}\left(\prod_{j=1}^{K_n}\frac{\lambda_{I_n}}{\lambda^{\mathbb Q}_{I_n}}\frac{f^{\mathbb Q}_{I_n}}{f_{I_n}}\cdot \frac{\alpha_{I_n,+}}{\alpha_{I_n,+}^{\mathbb Q}}{\rm e}^{-(\alpha_{I_n,+}-\alpha_{I_n,+}^{\mathbb Q})A_{I_n,j}^+}\frac{\alpha_{I_n,-}}{\alpha_{I_n,-}^{\mathbb Q}}\cdot{\rm e}^{-(\alpha_{I_n,-}-\alpha_{I_n,-}^{\mathbb Q})A_{I_n,j}^-}\cdot{{\rm e}^{-\omega\s B_{I_n,j}}}{b_{I_n}(-\omega\s)}\right)
\\
&\:\:\:\hspace{1.7cm}\times\left(\frac{q_{I_n}}{q^{\mathbb Q}_{I_n}}\frac{f^{\mathbb Q}_{I_n}}{f_{I_n}} \cdot \frac{\alpha_{I_n,+}}{\alpha_{I_n,+}^{\mathbb Q}}{\rm e}^{-(\alpha_{I_n,+}-\alpha_{I_n,+}^{\mathbb Q})A_{I_n,K_n+1}^+}\frac{\alpha_{I_n,-}}{\alpha_{I_n,-}^{\mathbb Q}}\cdot{\rm e}^{-(\alpha_{I_n,-}-\alpha_{I_n,-}^{\mathbb Q})A_{I_n,K_n+1}^-}\right);
\end{align*}
to understand this expression, recognize the effect of drawing the initial state, the $K_n$ claim arrivals (and the maxima in the corresponding intervals) before time $T_n$, and the event that $T_n$ occurs before a possible $(K_n+1)$-st arrival (and the maximum in the corresponding interval).
It requires some elementary (but rather tedious) algebra to check that, for any $i\in\{1,\ldots,d\}$, 
\[(\alpha_{i,+}-\alpha_{i,+}^{\mathbb Q})A_{i,j}^+ + (\alpha_{i,-}-\alpha_{i_n,-}^{\mathbb Q})A_{i,j}^-=-\omega\s
(A_{i,j}^+ -A_{i,j}^-);\]
in addition,
\[\alpha_{i,+} \alpha_{i,-} = 2\frac{f_i}{\sigma_i^2},\:\:\:\:\:\:\alpha^{\mathbb Q}_{i,+} \alpha^{\mathbb Q}_{i,-} = 2\frac{f^{\mathbb Q}_i}{\sigma_i^2},\:\:\:\:\:
\frac{\lambda_i}{\lambda_i^{\mathbb Q}}\cdot{b_i(-\omega\s)}=1,\:\:\:\:\:\: \frac{p_{i}}{p_{i}^{\mathbb Q}}\cdot
\frac{q_{i}}{q^{\mathbb Q}_{i}}=1.\]
It thus follows that, for $n\in\{1,\ldots,N-1\}$,
\[L_n = \exp\left(-\omega\s \sum_{j=1}^{K_n+1}(A_{I_n,j}^+-A_{I_n,j}^-)- \omega\s \sum_{j=1}^{K_n}B_{I_n,j} \right) = {\rm e}^{-\omega Y_n};\]
this is not surprising, given that the change-of-measure we set up corresponds to exponentially twisting the $Y_n$s (and ${\mathbb E}\,{\rm e}^{\omega\s Y_n}=1$). 

We now shift our attention to the contribution to $L$ due to $(T_{N-1},T_N].$
Similar to the computations performed above, we obtain the following expression for $L_N$. There are two scenarios. In the first place, $u$ can be reached by a claim arrival; say this happens due to the $\bar K_N$-th claim arrival in the interval $(T_{N-1},T_N].$ Then, with $Z_N(u):=X(\tau(u))-X(T_{N-1})$,
\[L_N = \frac{p_{I_N}}{p_{I_N}^{\mathbb Q}} \exp\left(-\omega\s \sum_{j=1}^{\bar K_N}(A_{I_N,j}^+-A_{I_N,j}^-)- \omega\s \sum_{j=1}^{\bar K_N}B_{I_N,j} \right) =\frac{p_{I_N}}{p_{I_N}^{\mathbb Q}}{\rm e}^{-\omega\s Z_N(u)} .\]
In the second place, $u$ can be reached in between two claim arrivals; say this happens between the $\bar K_N$-th and $(\bar K_N+1)$-st claim arrival in $(T_{N-1},T_N].$ Now,
\begin{align*}L_N &= \frac{p_{I_N}}{p_{I_N}^{\mathbb Q}} \,\gamma_{I_N}\exp\left(-\omega\s \sum_{j=1}^{\bar K_N}(A_{I_N,j}^+-A_{I_N,j}^-)- \omega\s \sum_{j=1}^{\bar K_N}B_{I_N,j} -\omega\s A^+_{I_N,\bar K_N+1}\right)\\&=\frac{p_{I_N}}{p_{I_N}^{\mathbb Q}} \,\gamma_{I_N}{\rm e}^{-\omega\s Z_N(u)},\end{align*}
with
\[\gamma_i:= \frac{\lambda_{i}}{\lambda^{\mathbb Q}_{i}}\frac{f^{\mathbb Q}_{i}}{f_{i}}\cdot \frac{\alpha_{i,+}}{\alpha_{i,+}^{\mathbb Q}}=\frac{1}{b_i(-\omega\s)}\cdot \frac{\alpha_{i,-}^{\mathbb Q}}{\alpha_{i,-}}
.\]
Now observe that, by definition of $\tau(u)$, 
\[\sum_{n=1}^{N-1} Y_n+ Z_N(u) = X(T_{N-1})+ \big(X(\tau(u))-X(T_{N-1})\big)= 
X(\tau(u)) \geqslant u.\] Define
\[\Omega := \max_{i\in\{1,\ldots,d\}}\left(\frac{p_{i}}{p_{i}^{\mathbb Q}}\max\{\gamma_i,1\}\right) =
\max_{i\in\{1,\ldots,d\}}\left(\frac{q-\varphi_{i}(-\omega\s)}{q}\max\{\gamma_i,1\}\right).\]
The above yields that $L\leqslant \Omega \,{\rm e}^{-\omega\s u} $ almost surely.  In particular, we have derived the following Lundberg-type inequality for the resampling model.
\begin{theorem} \label{THUB} For any $u\geqslant 0$,
\[\pi(u) \leqslant \Omega\, {\rm e}^{-\omega\s u} .\]
\end{theorem}
\begin{remark}{\em Related Lundberg-type inequalities have appeared; cf.\ Remark \ref{RA}. 
We refer to \cite[Corollary 3.6]{AA} for a result covering the case that the L\'evy processes are of compound Poisson type.}
\end{remark}

Because of the almost sure upper  bound on $L$, we have also proven an optimality property of our importance sampling algorithm, namely that it has bounded relative error. This claim follows directly from the observation that ${\mathbb V}{\rm ar}_{\mathbb Q}(L)\leqslant {\mathbb E}_{\mathbb Q}(L^2)\leqslant \Omega^2 {\rm e}^{-2\omega\s u}.$ 
The definition of `bounded relative error' is provided in e.g. \cite[Section VI.1]{AG}. If an estimation procedure has bounded relative error, then  this effectively entails that the number of simulation runs needed to obtain an estimate with a given precision (e.g.\ 10\%), is bounded in $u$.
\begin{theorem} \label{BRE} If $\pi(u)\,{\rm e}^{\omega\s u}\to A>0$ as $u\to\infty$, then the procedure given by Algorithm $\ref{ALG1}$ has bounded relative error. 
\end{theorem}

\section{Examples, numerics}\label{EN}
In this section we focus on examples corresponding to the case $d=2.$ In the first subsection we provide explicit computations pertaining to the quantities that play a role in our exact asymptotics, whereas in the second subsection we present  illustrative numerical examples, which in particular  quantify the potential risk due to ignoring the resampling. 

\subsection{Explicit expression for two-dimensional case}
We now point out how to compute various quantities needed to evaluate the exact asymptotics of Thm.\ \ref{THCL}. We let, as justified before,  both L\'evy processes $X_1(\cdot)$ and $X_2(\cdot)$ be sums of Brownian motions and compound Poisson processes. 
\begin{itemize}
\item[$\circ$]
As a first step, we have to find the two solutions (for $\theta$) to the equation 
\[\frac{1}{q} = \frac{p_1}{q-\varphi_1(\alpha)+\theta}+\frac{p_2}{q-\varphi_2(\alpha)+\theta},\]
where we assume that $p_1\varphi_1'(0)+p_2\varphi_2'(0)>0.$
After some elementary algebra, one finds that this equation can be rewritten as 
\[\theta^2 -\theta(\varphi_1(\alpha)+ \varphi_2(\alpha)-q) + \varphi_1(\alpha)\,\varphi_2(\alpha)-\varphi_1(\alpha) \,p_1-
\varphi_2(\alpha)\,p_2=0.\]
With $\Delta(\alpha)$ defined as
\[(\varphi_1(\alpha)-\varphi_2(\alpha))^2 -2q(\varphi_1(\alpha)+\varphi_2(\alpha))+q^2+4q(\varphi_1(\alpha)\,p_1+\varphi_2(\alpha)\,p_2),\]
we thus obtain
\[\theta_{k}(\alpha) =\tfrac{1}{2}\big(\varphi_1(\alpha)+ \varphi_2(\alpha)-q\big)\pm \tfrac{1}{2}\sqrt{\Delta(\alpha)} \]
(where we let the expression with the minus-sign correspond to $k=1$, and the expression with the plus-sign to $k=2$, so that $\bar\theta(\alpha)=\theta_2(\alpha)$).
As we saw before, $\theta_1(\alpha)\in(\bar\varphi_1(\alpha), \bar\varphi_2(\alpha))$ and $\theta_2(\alpha)>\bar\varphi_2(\alpha).$
\item[$\circ$]
We now compute the matrices $S(-\omega\s)$ and $S^{-1}(-\omega\s)$. 
From Remark \ref{R2} we know that the $j$-th component of the $k$-th eigenvector of $M(\alpha)$ is given by
\[S_{jk}(\alpha) =\frac{1}{q-\varphi_j(\alpha)+\theta_k(\alpha)}=\frac{2}{q-\varphi_j(\alpha)+\varphi_{3-j}(\alpha)+(-1)^k\sqrt{\Delta(\alpha)}},\]
for $j=1,2$ and $k=1,2$. In addition, with $D(\alpha):=S_{11}(\alpha)S_{22}(\alpha)-S_{21}(\alpha)S_{12}(\alpha),$
\[S^{-1}(\alpha)=\frac{1}{D(\alpha)}\left(\begin{array}{rr}S_{22}(\alpha)&-S_{12}(\alpha)\\-S_{21}(\alpha)&S_{11}(\alpha)\end{array}\right).\]
\item[$\circ$]
The decay rate $-\omega\s$ is the negative solution to the equation $\theta_2(\alpha)=0$. By squaring 
$\varphi_1(\alpha)+ \varphi_2(\alpha)-q=\sqrt{\Delta(\alpha)}$, this is a
negative solution to
\begin{equation}\label{fp}\varphi_1(\alpha)\,\varphi_2(\alpha)-q\,\varphi_1(\alpha)\,p_1-q\,\varphi_2(\alpha)\,p_2=0\end{equation}
(but because of the squaring we have to verify whether we found an admissible root). 

\item[$\circ$]
We obtain, using the notation in Eqn.\ (\ref{AA}),
\[{\bs u}=\left(\begin{array}{r}S_{12}(-\omega\s)\\S_{22}(-\omega\s)\end{array}\right),\:\:\:{\bs v}=\frac{1}{D(-\omega\s)}\left(\begin{array}{r}-S_{21}(-\omega\s)\\S_{11}(-\omega\s)\end{array}\right).\]

\item[$\circ$]
To evaluate the constant $A$ featuring in Thm.\  \ref{THCL}, we are left with computing the vector ${\bs\ell}$. 
This we do by computing the generator matrix $\Lambda$ pertaining to the process $(J(\eta(v)))_{v\geqslant 0}$, as was introduced in Remark \ref{R1}, applying the results of \cite{DIKM}; note that in principle also \cite{IM} can be used, as we are dealing with a time-reversible process $X(\cdot)$. One approach is to rely on the matrix integral equation discussed in \cite[Section 4.1]{DIKM}, but we here apply the more explicit characterization of \cite[Thm.\ 1]{DIKM}, as follows. 

\noindent
Observe that in our case $d^-=0$, as there are no states corresponding to a decreasing subordinator. 
For our two-dimensional setting this means that the results of \cite{DIKM} entail that we can write the transition rate matrix $\Lambda$ can be written as $-V\Gamma V^{-1}$, where $\Gamma$ is the diagonal matrix with the non-negative zeroes of ${\rm det}(M(\alpha))$ on the diagonal, and the columns of $V$ consist of the corresponding right eigenvalues; as a consequence of $d^-=0$ we have, 
in the terminology of \cite{DIKM},  that $V=V^+$. 

\noindent The next step is to consider the non-negative roots of  ${\rm det}(M(\alpha))$. First note that the matrix $M(\alpha)$ is given by
\[M(\alpha) =\left(\begin{array}{cc}-qp_2+\varphi_1(\alpha)&qp_2\\qp_1&-qp_1+\varphi_2(\alpha)\end{array}\right).\]
The corresponding determinant can be written as, with $m_1(\alpha):=\varphi_1(\alpha)\,\varphi_2(\alpha)$ and $m_2(\alpha):=
p_1\,\varphi_1(\alpha)+p_2\,\varphi_2(\alpha)$,
\[m(\alpha):={\rm det}(M(\alpha))=m_1(\alpha)-q\,m_2(\alpha)=0.\]
Obviously, $0$ is a root of this equation. In addition, e.g.\ \cite[Corollary 5]{IM} states that there is precisely one positive root, which we call $\alpha\s$. As a consequence of the facts that (i) at least one of the $\varphi_i(\alpha)$ is positive for all $\alpha>0$, (ii) $m_1'(0)=0$ and $m_2'(0)=\kappa>0$, and (iii) $m_2(\alpha)>0$ for all $\alpha>0$ (due to the convexity of $m_2(\cdot)$), and (iv) $m_1(\alpha)/m_2(\alpha)\to\infty$ as $\alpha\to\infty$,  it follows that necessarily
$\varphi_1(\alpha\s)>0$ and $\varphi_2(\alpha\s)>0$.

\noindent
We thus obtain
\[\Gamma=\left(\begin{array}{cc}0&0\\0&\alpha\s\end{array}\right),\:\:\:
V=\left(\begin{array}{cc}1&qp_2\\1&qp_2-\varphi_1(\alpha\s)\end{array}\right),\]
so that
\begin{align*}\Lambda&=\frac{1}{\varphi_1(\alpha\s)}
\left(\begin{array}{cc}1&qp_2\\1&qp_2-\varphi_1(\alpha\s)\end{array}\right)
\left(\begin{array}{cc}0&0\\0&\alpha\s\end{array}\right)
\left(\begin{array}{cc}qp_2-\varphi_1(\alpha\s)&-qp_2\\-1&1\end{array}\right)
 \\
&= \frac{\alpha\s\,q}{\varphi_1(\alpha\s)\,\varphi_2(\alpha\s)}\left(\begin{array}{rr}-p_2\,\varphi_2(\alpha\s)&p_2\,\varphi_2(\alpha\s)\\p_1\,\varphi_1(\alpha\s)&-p_1\,\varphi_1(\alpha\s)\end{array}\right)
;\end{align*}
here we have used  that (\ref{fp}) implies that $\varphi_2(\alpha\s)\,(\varphi_1(\alpha\s) -qp_2)=\varphi_1(\alpha\s)\,qp_1$. Applying this expression for $\Lambda$, elementary computations thus lead to
\[\bar{\bs \pi} = \frac{1}{p_1\,\varphi_1(\alpha\s)+p_2\,\varphi_2(\alpha\s)}\left(\begin{array}{r}
p_1\,\varphi_1(\alpha\s)\\p_2\,\varphi_2(\alpha\s)\end{array}\right).\]
(Alternatively, one could use the relation, with ${\bs 1}:=(1,0)^\top$ denoting the first unit vector, \[
\bar{\bs \pi} = \frac{1}{v} \cdot {\bs 1}^\top V^{-1},\]
with $v$ denoting the normalizing constant ${\bs 1}^\top V^{-1}{\bs e}.$)
\end{itemize}

Now consider more specifically the case  that (for $i=1,2$) $X_i(\cdot)$ is the sum of a Brownian motion and a compound Poisson process, so that we can write
\[\varphi_i(\alpha)= \tfrac{1}{2}\sigma_i^2\,\alpha^2 +r_i\,\alpha +\lambda_i\,(b_i(\alpha)-1);\]
here the variance coefficients $\sigma_i^2$ are positive, and in our insurance context typically the premium rates $r_i$ as well.
We have that the (negative of the) asymptotic drift $\kappa$, as was defined before, equals $p_1(r_1+\lambda_1\,b_1'(0))+ p_2(r_2+\lambda_2\,b_2'(0))$, which we have assumed to be positive.

Specializing to the case of exponentially distributed claims (such that $b_i(\alpha)=\mu_i/(\mu_i+\alpha)$, for some $\mu_i>0$), we have 
\[m(\alpha)=\prod_{i=1}^2\left(\tfrac{1}{2}\sigma_i^2\,\alpha^2 -r_i\,\alpha -\frac{\lambda_i\alpha}{\mu_i+\alpha}\right)-
q\sum_{i=1}^2\left(\tfrac{1}{2}\sigma_i^2\,\alpha^2 -r_i\,\alpha -\frac{\lambda_i\alpha}{\mu_i+\alpha}\right)p_i=0,\]
which is (after some rewriting) a polynomial equation of degree 6 that can be solved by standard software, so as to obtain $-\omega\s$ and $\alpha\s$. 

\subsection{Numerical example} 

For the numerical results we have used a setup that is as much as possible in line with the one considered in  \cite{ASMS}. 
\begin{itemize}
\item[$\circ$]
The environmental process has stationary distribution ${\bs p}=\left(\frac{2}{3},\frac{1}{3}\right)$ and the intensity $q$ is of the form $3\cdot 4^i$, for $i\in\{-2,-1,0,1,2\}$. 
\item[$\circ$]
We let the premium rate, the variance coefficient of the Brownian terms,  and the claim sizes be environment-independent: the premium rates are ${\bs r}=(1,1)$, the Brownian motions are characterized by ${\bs \sigma}=(1,1)$, and the claim sizes are exponentially distributed claims with parameter ${\bs \mu}=(1.1, 1.1)$. 
\item[$\circ$]
The intensities ${\bs \lambda}$  of the claims sizes are chosen, again following Asmussen \cite{ASMS}, such that $\lambda_1=\rho/2$ and $\lambda_2=2\rho$ where $\rho$ denotes the average amount of claim per unit time, i.e., \[\rho=-\sum_{i=1}^2p_i\lambda_ib'_i(0)=\sum_{i=1}^2 p_i\frac{\lambda_i}{\mu_i}.\] The value of $\rho$ is fixed at 0.9, so that we have ${\bs \lambda}=(0.45,1.8)$. 
\end{itemize}

In Table \ref{table_env} we present the corresponding numerical output. 
The column `Exact' is the value of $\pi(u)$ determined by an importance-sampling based computation;  the algorithm presented in Section \ref{UNI} (which has bounded relative error) has been used. The fact that, per parameter setting, we have used as many as $200\,000$ runs guarantees estimates with a high precision. The next two columns present the exact asymptotics $A\,{\rm e}^{-\omega\s u}$ of Thm.\ \ref{THCL} and the upper bound $\Omega\,{\rm e}^{-\omega\s u}$ of Thm.\ \ref{THUB}, respectively. The last column provides the values that one would get if the resampling were ignored. Then the input process is a (non-modulated) L\'evy process with arrival rate $\bar{\lambda}=p_1\lambda_1+p_2\lambda_2$, so that $\pi(u)$ can be evaluated using results presented in e.g.\ \cite{ASM,DM}. Alternatively, these values can be found by taking  $q$ sufficiently large in the modulated model.

\begin{table}[h]
\centering
\begin{tabular}{|l|cccc|}
\hline
$u=175$ &Exact & Thm.\ \ref{THCL} & Thm.\ \ref{THUB} & No modulation \\
 \hline
$q=3\cdot 4^{-2}=0.1875$ & $9.21\cdot 10^{-3}$ & $9.21\cdot 10^{-3}$  & $1.12\cdot 10^{-2}$& $6.26\cdot 10^{-6}$\\
$q=3\cdot 4^{-1}=0.75$& $1.90\cdot 10^{-4}$ & $1.89\cdot 10^{-4}$ & $2.11\cdot 10^{-4}$ & $6.26\cdot 10^{-6}$\\
$q=3\cdot 4^{0}=3$ & $1.86\cdot 10^{-5}$ & $1.86\cdot 10^{-5}$ & $1.98\cdot 10^{-5}$ & $6.26\cdot 10^{-6}$\\
$q=3\cdot 4^{1}=12$ & $8.36\cdot 10^{-6}$ & $8.36\cdot 10^{-6}$ & $8.80\cdot 10^{-6}$& $6.26\cdot 10^{-6}$ \\
$q=3\cdot 4^{2}=48$ & $6.72\cdot 10^{-6 }$& $6.72\cdot 10^{-6}$ & $7.05\cdot 10^{-6}$& $6.26\cdot 10^{-6}$ \\
\hline
\end{tabular}
\caption{\label{table_env}Numerical results varying the speed $q$ of the background process.}
\end{table}

The conclusions from the above table are the following. 
\begin{itemize}
\item[$\circ$]
In the first place, a comparison between the columns `Exact' and `No modulation' reveals that by ignoring resampling one potentially substantially underestimate the risk, in particular when the timescale of resampling is slow relative to the timescale corresponding to (the jumps of) the L\'evy processes $X_1(\cdot)$ and $X_2(\cdot)$; this corresponds to the regime of small $q$. In the regime that $q$ is relatively large, we observe that the resampling is apparently so frequent that the individual L\'evy processes can be safely replaced by their time-average counterpart. 
\item[$\circ$]In addition, it is observed that the approximation based on Thm.\ \ref{THCL} is nearly exact. The upper bound based on Thm.\ \ref{THUB} is typically rather tight (in the table the relative error is between $5\%$ and $20\%$), particularly when $q$ is relatively large. 
\item[$\circ$] The numerical output also confirms the (intuitively clear) property  that adding the Brownian component (with ${\bs \sigma}=(1,1)$) leads to a higher ruin probability; this follows by comparing our output with that presented  in \cite{ASMS}, in which the claim process does not contain a Brownian component.
\end{itemize}

In Table \ref{table_u} we fix the environmental intensity  $q$ at $\frac{3}{4}$ and vary the value of $u$. 
Again, it is seen that ignoring the resampling may lead to a significant underestimation of the ruin probability. The other conclusions are similar to those corresponding to Table \ref{table_env}. In this parameter setting (i)~$A/\Omega$ is approximately $0.9$ and (ii)~the approximation based on Thm.~\ref{THCL} is near-exact, thus entailing that the upper bound of Thm.~\ref{THUB} is about 10\% off. 

\begin{table}[h]
\centering
\begin{tabular}{|l|cccc|}
\hline
$q=\frac{3}{4}$ &Exact & Thm.\ \ref{THCL} & Thm.\ \ref{THUB} & No modulation \\
 \hline
$u=175$ & $1.90\cdot 10^{-4}$ &$1.89\cdot 10^{-4}$ & $2.11\cdot 10^{-4}$ & $6.26\cdot 10^{-6}$ \\
$u=162.5$ & $3.48\cdot 10^{-4}$ & $3.47\cdot 10^{-4}$ &  $3.87\cdot 10^{-4}$ & $1.47\cdot 10^{-5}$\\
$u=150$ & $6.38\cdot 10^{-4}$ & $6.37\cdot 10^{-4}$ & $7.10\cdot 10^{-4}$ & $3.44\cdot 10^{-5}$ \\
$u=137.5$ & $1.17\cdot 10^{-3}$ & $1.17\cdot 10^{-3}$ & $1.30\cdot 10^{-3}$ & $8.07\cdot 10^{-5}$ \\
$u=125$ & $2.15\cdot 10^{-3}$ & $2.14\cdot 10^{-3}$& $2.39\cdot 10^{-3}$ &  $1.89\cdot 10^{-4}$\\
\hline
\end{tabular}
\caption{\label{table_u}Numerical results varying the initial reserve $u$.}
\end{table}

\section{Discussion and concluding remarks}

This paper addresses the evaluation of ruin probabilities for a model in which the underlying dynamics are periodically resampled. We have generalized two celebrated results from the risk theory literature: we identify the exact tail asymptotics (thus extending `Cram\'er-Lundberg'), and derived a uniform upper bound (thus extending Lundberg's inequality). In our proof of  the uniform upper bound, we developed an importance-sampling-based efficient simulation algorithm which is proven to have bounded relative error. Numerical experiments showed that neglecting the parameter uncertainty typically leads to a significant underestimation of the ruin probability.

Various extensions can be thought of; we mention three directions for future research. Where this paper focuses on the probability of ultimate ruin, a first obvious extension would relate to the {\it finite-horizon} ruin probability (i.e., ruin before some $T>0$). In the importance sampling procedure one would anticipate that one should distinguish between the case in which ultimate ruin corresponds (with high probability) with a ruin time smaller than $T$ (such that the change of measure of Section~\ref{UNI} can be used), and the case in which ultimate ruin corresponds (with high probability) with a ruin time larger than $T$ (such that the claim arrival process should be `twisted' more strongly). It is expected that the same dichotomy appears in the exact asymptotics and the uniform upper bound. 

In the second place, one could aim at relaxing the exponentiality assumptions. More concretely, one could consider the model in which there is resampling at phase-type \cite[Section IX.1]{AA}  distributed times. Likewise, an interesting extension concerns generalizing the results for the case that the processes $X_i(\cdot)$ are compound Poisson (with, for each $i$, exponentially distributed claim interarrival times) to their counterparts in which  the claim interarrival times are of phase-type. Though notationally rather involved, conceptually such extensions are relatively straightforward; see e.g.\ \cite{KM} for such computations in a related model. In addition, one could consider the specific case of Gamma claims; cf.\ \cite{CSZ}.

In the third place one could consider a model involving multiple business lines such that the individual claim arrival process react to a {\it common} environmental process. In this setup there is correlation between the ruin events; one could for instance aim at computing the probability of ruin of (minimally) one of the business lines, or alternatively the probability of ruin of all of them. A relevant related problem concerns the allocation of a firm's capital to the individual business lines; cf.\ the model considered in \cite{DMSW}. 

\bibliographystyle{plain}

\end{document}